\documentclass[11pt]{article}
\usepackage{geometry}                
\usepackage{graphicx}
\usepackage{amssymb}
\usepackage{epstopdf}
\usepackage{amsmath,amsthm,amscd,amssymb}
\usepackage{latexsym}
\numberwithin{equation}{section}

\theoremstyle{plain}
\newtheorem{theorem}{Theorem}[section]
\newtheorem{lemma}[theorem]{Lemma}

\theoremstyle{definition}
\newtheorem{definition}[theorem]{Definition}

\theoremstyle{remark}

\newtheorem{case[theorem]}{Case}

\def\norm#1.#2.{\lVert#1\rVert_{#2}}

\title{Ubiquity of simplices in subsets of vector spaces over finite
fields}
\author{Derrick Hart and Alex Iosevich}

\begin{document}

\maketitle

\begin{abstract} We prove that a sufficiently large subset of the $d$-dimensional 
vector space over a finite field with $q$ elements, $ {\Bbb F}_q^d$, contains a copy
 of every $k$-simplex. Fourier analytic methods, Kloosterman sums, and bootstrapping play an important role.
\end{abstract}

\tableofcontents

\section{Introduction}

\vskip.125in
Many problems in combinatorial geometry ask, in one form or another,
whether a certain structure must be present in a set of sufficiently
large size. Perhaps the most celebrated result of this type is
Szemeredi's theorem (\cite{SO75}) which says that if a subset of the
integers has positive density, then it contains an arbitrary large
arithmetic progression. The conclusion has recently been extended to
the subsets of prime numbers by Green and Tao (\cite{GT07}). In
Euclidean space, a result due to Katznelson and Weiss (\cite{FKB90})
says that a subset of Euclidean space of positive Lebesgue upper
density contains every sufficiently large distance. A subsequent
result by Bourgain (\cite{BK86}), says that a subset of $\Bbb R^k$ of
positive Lebesgue upper density contains an isometric copy of all
large dilates of a set of $k$ points spanning a $(k-1)$-dimensional
hyperplane. Ergodic theory has been used to show that positive upper
density implies that the set contains a copy of a sufficiently large
dilate of every convex polygon with finitely many sides. See, for
example, a recent survey by Bryna Kra (\cite{K06}).

Let $\Bbb F_q^d$ be a $d$-dimensional vector space over a finite
field $\Bbb F_q$ of odd characteristic.  A plausible analogy to Bourgain's result (\cite
{BK86}) in this context would be to consider whether a subset of
positive density contains a isometric copy of a set of $k$ points
spanning a $(k-1)$-dimensional hyperplane. It turns out however, that
the positive density condition is much too strong in the context of vector spaces over finite fields and the same
conclusion follows from a much weaker assumption on the size of the
underlying set.

\begin{definition}
Let a $k$-simplex be a set of $k+1$ points in general position, which
means that
no $n + 1$ of these points, $n \leq k$, lie in a $(n - 1)$-dimensional
sub-space of ${\Bbb F}_q^d$.
\end{definition}

\begin{definition} We say that a linear transformation $T$ on ${\Bbb
F}_q^d$ is an isometry if
$$ ||Tx||=||x||,$$ where
$$ ||x||=x_1^2+x_2^2+\dots+x_d^2,$$ an element of ${\Bbb F}_q$.
\end{definition}

The question we ask in this paper is how large does $E \subset {\Bbb
F}_q^d$ need to be in order to be sure that it contains a copy of
every $k$-simplex.
Our main result is the following.
\begin{theorem} \label{main} Let $E \subset {\Bbb F}_q^d$, $d>{k+1 \choose 2}$, such that
$|E| \ge C q^{\tfrac{k}{k+1}d}q^{\tfrac{k}{2}}$ with a
sufficiently large constant $C>0$.  Then $E$ contains an isometric
copy of every $k$-simplex.
\end{theorem}

Note that we obtain non-trivial results only when $k<<\sqrt{d}$.
Nevertheless, in that range we are able to dip considerably below the
positive density condition on the underlying set $E$.

The method of proof relies on the fact that orthogonal
transformations on ${\Bbb F}_q^d$ are isometries.
A "distance representation" of a simplex is then used to reduce
Theorem \ref{main} to an appropriate weighted incidence theorem for
spheres and points.  Weil's estimate (\cite{We48}) for classical
Kloosterman sums is used to control the size of the Fourier transform
of spheres of non-zero radius. The key idea in the proof is to show
at each step of an inductive argument that a collection of distances
among vertices of a given simplex can not only be realized, but
actually occur a "statistically correct" number of times.
\vskip.125in

\section{Preliminaries and Definitions}
Let $\Bbb F _q^d$ be the $d$-dimensional vector space over the finite
field ${\Bbb F}_q$. The Fourier transform of a function
$$f: \Bbb F_q^d \rightarrow \Bbb C$$ is given by
\begin{equation*}
\widehat{f}(m) := q^{-d} \sum_{x \in \Bbb F_q^d} f(x) \chi(-x \cdot m),
\end{equation*}
where $\chi$ is an additive character on $\Bbb F_q$.

The orthogonality property of the Fourier Transform says that
$$ q^{-d}\sum_{x \in \Bbb F_q^d} \chi(-x \cdot m)=1$$ for $m=(0,
\dots, 0)$ and $0$
otherwise yields many standard properties of the Fourier Transform.

We summarize some of the properties of the Fourier Transform as follows.

\begin{lemma}[The Fourier Transform]
Let
$$f,g:\Bbb F_q^d \rightarrow \Bbb C.$$

\begin{align*}
& \hat f (0, \dots, 0) =q^{-d} \sum_{x \in \Bbb F_q^d} f(x),
\\
& (\text{Plancherel})\ q^{-d} \sum_{x \in \Bbb F_q^d} f(x) \overline{g
(x)} =\sum_{m \in \Bbb F_q^d} \hat{f}(m) \overline{\hat{g}(m)},
\\
& (\text{Inversion})\ \ \  f(x) =\sum_{m \in \Bbb F_q^d} \hat{f}(m)
\chi(x \cdot m).
\end{align*}
\end{lemma}

\subsection{Notation} Throughout the paper $X \lesssim Y$ means that
there exists $C>0$ such that $X\leq CY$, $X \gtrsim Y$ means $Y
\lesssim X$, and $X \approx Y$ if both $X\lesssim Y$ and
$X\gtrsim Y$.  Along the same lines, $X \ll Y $ means that $\tfrac{X}
{Y}\rightarrow 0$, as $q \to \infty$ $X \gg Y$ means $Y\ll X$, and $X
\sim Y$ if $\tfrac{X}{Y}\rightarrow 1$ as $q \rightarrow \infty$.
\vskip.125in

\section{Proof of the main result}
Even though a finite field with $q$ elements, ${\Bbb F}_q$, is not a
metric space, we define the "distance" between two points $x$ and $y$
in ${\Bbb F}_q^d$ by the formula
$$ ||x-y||={(x_1-y_1)}^2+{(x_2-y_2)}^2+\dots+{(x_d-y_d)}^2.$$
The same notion of "distance" was used by Bourgain, Katz and Tao(\cite
{BKT04}), and Iosevich and Rudnev(\cite{IR07}) in their study of the
Erd\H os distance problem in vector spaces over finite fields. As we
noted above, a geometric justification of this notion of distance is
that an orthogonal transformation on ${\Bbb F}_q^d$, a matrix $O$
such that $O^t \cdot O=I$, preserves this notion of a "distance".
Represent a $k$-simplex in a subset $E \subset \Bbb F_q^d$ on $k+1$
points recursively by setting
$$ {\cal T}_{l_k}=\{(x_0,\dots,x_{k-1},x_k) \in {\cal T}_{l_{k-1}}
\times E: ||x_0-x_k||=t_{1,k}, ||x_1-x_k||=t_{2,k},\dots,||x_{k-1}-x_k||=t_{k,k}\},$$
for $l_k=l_{k-1}\cup_\{t_{1,k},\dots t_{k,k}\},t_{i,j}\in \mathbb F_q^*$ where
$$ {\cal T}_{l_1}=\{(x_0,x_1) \in E^2: ||x_0-x_1||=t_{1,1}\}.$$
This representation does not, in general, always embody a simplex as $
{\cal T}_{l_k}^k$ is not guaranteed to be in general position.
However, as we show below, "legitimate" $k$-simplices are equivalent
up to an orthogonal transformation.

\begin{theorem} \label{simplexlemma} Let $E \subset {\Bbb F}_q^d$, $d>
{k+1 \choose 2}$, such that
$|E| \ge Cq^{\tfrac{k}{k+1}d}q^{\tfrac{k}{2}}$, with a
sufficiently large constant $C$.  Then for every side length set $l_k$, $l_k\in (\Bbb F_q^*)^{k+1 \choose 2}$ we have that $|{\cal T}_{l_k}|>0$.  Furthermore,
$$ |{\cal T}_{l_k}| \sim |E|^{k+1}q^{-{k+1 \choose 2}}.$$
\end{theorem}

Using this theorem we recover the main result of the paper using the
following linear algebraic observation.

\begin{lemma}\label{uptocrap} Let $P$ be a simplex with vertices
$V_0, V_1, \dots, V_k$, $V_j \in {\Bbb F}_q^d$. Let $P'$ be another
simplex with vertices $V'_0, V'_1, \dots, V'_k$. Suppose that
\begin{equation} \label{equalnorm} ||V_i-V_j||=||V'_i-V'_j|| \end{equation}
for all $i,j$.
Then there exists an orthogonal, affine transformation $O$ on ${\Bbb F}_q^d$ such
that $O(P)=P'$.
\end{lemma}

\vskip.125in

\subsection{Proof of Theorem \ref{simplexlemma}-the main result
reformulated in terms of "distances"}

The proof proceeds by induction. The first step is the case $k=2$.
For a set $E$ we define the characteristic or indicator function to
be $E(x)$.  Now define the sphere of radius $t_{1,1} \in \Bbb F_q^*$
to be
$$S_{t_{1,1}}=\{x \in {\Bbb F}_q^d: ||x||=t_{1,1}\},$$ then
\begin{align*}
|{\cal T}_{l_1}|&=|\{(x_0,x_1) \in E \times E: ||x_0-x_1||=t_{1,1}\}|
\\
&=\sum_{x_0,x_1} E(x_0)E(x_1)S_{t_{1,1}}(x_0-x_1).
\end{align*}

In order to obtain information from this quantity the behavior of
incidences of spheres and points in $E$ will be critical.  The
following classical fact, whose proof will be given in a subsequent
section, states that the sphere has optimal Fourier decay away
from the origin.

\begin{lemma} \label{sphere} Let $S_{t}$, $t \in \Bbb F_q^*
$ be defined as above. If $m \not=(0, \dots, 0)$ then
$$ |\widehat{S}_{t}(m)| \lesssim q^{-\frac{d+1}{2}}, $$ and
$$ \widehat{S}_{t}(0, \dots, 0)=q^{-d} |S_{t}| \approx q^{-1}.$$ \end{lemma}
Applying Fourier inversion to the sphere,
\begin{align*}
|{\cal T}_{l_1}|=&q^{2d}\sum_{x_0,x_1} E(x_0)E(x_1)\sum_m \widehat{S}_{t_{1,1}}(m)\chi(m\cdot(x_0-x_1))
\\
&=q^{2d} \sum_m {|\widehat{E}(m)|}^2 \widehat{S}_{t_{1,1}}(m)
\\
&={|E|}^2 \cdot q^{-d} \cdot |S_{t_{1,1}}|+q^{2d} \sum_{m \not=(0,
\dots, 0)} {|\widehat{E}(m)|}^2 \widehat{S}_{t_{1,1}}(m)
\\
&=M+R.\,
\end{align*}
By Lemma \ref{sphere},
$$ M \approx \frac{{|E|}^2}{q},$$ and using Lemma \ref{sphere} once
again,
\begin{align*}
|R| &\lesssim q^{2d} \cdot q^{-\frac{d+1}{2}} \cdot \sum_m {|\widehat
{E}(m)|}^2
\\
&=q^{\frac{d-1}{2}} \cdot |E|,
\end{align*}
which is smaller than $M$ if $|E| \ge Cq^{\frac{d+1}{2}}$ with a
sufficiently large constant $C$ and thus
${\cal T}_{l_1}$ is non-empty. Moreover, if $|E| \gg q^{\frac{d+1}
{2}}$, we get the "statistically expected" number of distances,
$$|{\cal T}_{l_1}| \sim \frac{{|E|}^2}{q}.$$

Assuming the $(k-1)$st case, we count the number of $k$-simplices
in $E$ as an extension of the $(k-1)$-simplices in $E$.
$$|{\cal T}_{l_{k}}|= \sum_{x_0,\dots,x_k} {\cal T}_{l_{k-1}}(x_0,\dots,x_{k-1})E(x_k)S_{t_{1,k}}(x_0-x_k)\dots S_{t_{k,k}}(x_{k-1}-x_k).$$
By Fourier inversion, the expression equals
$$ \sum_{x_0,\dots,x_k} \sum_{m_0,\dots,m_{k-1}} \prod_{i=1}^{k}\chi
((x_{i-1}-x_k) \cdot m_{i-1})\widehat{S}_{t_{i,k}}(m_{i-1}) {\cal T}_{l_{k-1}}(x_0,\dots,x_{k-1})E(x_k)$$
$$=q^{(k+1)d} \sum_{m_0,\dots,m_{k-1}} \widehat{{\cal T}}_{l_{k-1}}
(-m_0,\dots,-m_{k-1}) \widehat{E}(m_0+\dots+m_{k-1})  \widehat{S}_
{t_{1,k}}(m_0)\dots\widehat{S}_{t_{k,k}}(m_{k-1}),$$ where the
Fourier transform of ${\cal T}_{l_{k-1}}$ is actually the
Fourier transform on ${\Bbb F}_q^d \times\dots\times {\Bbb F}_q^d$,
$k$ times.

Extracting the zero term and breaking the remaining sum into pieces
on which we may apply Lemma \ref{sphere}, this expression equals
$$q^{(k+1)d} \cdot q^{-(k+1)d} \cdot |{\cal T}_{l_{k-1}}| \cdot |E|
\cdot |S_{t_{1,k}}| \cdot q^{-d}\cdot\dots \cdot |S_{t_{k,k}}|
\cdot q^{-d}$$
$$+q^{(k+1)d} \sum_{\substack{\mathcal{I} \cup \mathcal{I}'= \{ 0, \dots ,k-1 \} \\ m_i=0\ (i\in \mathcal{I})\\ m_i\neq 0\ (i \notin \mathcal{I})}} 
\widehat{{\cal T}}_{l_{k-1}}(-m_0,\dots,-m_{k-1}) \widehat{E}(m_0+\dots+m_{k-1})\widehat{S}_{t_{1,k}}(m_0)\dots\widehat{S}_{t_{k,k}}(m_{k-1})$$
$$=M+R,$$ where the sum defining $R$ runs over all the partitions of $\{0,\dots,k-1\}$ with the case $\mathcal{I}'=\emptyset$
extracted and used as the main term $M$ above.

By Lemma \ref{sphere} and the induction hypothesis,
$$ M \sim {|E|}^{k+1}q^{-{k+1 \choose 2}}.$$

By Lemma \ref{sphere} we have that
$$|R|\lesssim q^{(k+1)d}\sum_{\substack{\mathcal{I} \cup \mathcal{I}'= \{ 0, \dots ,k-1 \} \\ m_i=0\ (i\in \mathcal{I})\\ m_i\neq 0\ (i \notin \mathcal{I})}} 
q^{-|\mathcal{I}'|(d+1)/2-|\mathcal{I}|}
|\widehat{{\cal T}}_{l_{k-1}}(-m_0,\dots,-m_{k-1})| |\widehat{E}(m_0+\dots+m_{k-1})|.$$

Then for each term in the sum corresponding to a partition $\mathcal{I} \cup \mathcal{I}'$ we apply Cauchy-Schwarz,
$$\sum_{\substack{m_i=0\ (i\in \mathcal{I})\\ m_i\neq 0\ (i \notin \mathcal{I})}} 
|\widehat{{\cal T}}_{l_{k-1}}(-m_0,\dots,-m_{k-1})| |\widehat{E}(m_0+\dots+m_{k-1})|
\lesssim A^{1/2}B^{1/2}.$$
Applying Plancherel and the induction hypothesis,
$$A\leq\sum_{m_0,\dots,m_{k-1}} |\widehat{{\cal T}}_{l_{k-1}}(-m_0,
\dots,-m_{k-1})|^2 =q^{-kd}|{\cal T}_{l_{k-1}}|\sim q^{-kd}q^{-{k \choose 2}} {|E|}^{k}.$$
Now
$$B=\sum_{m_i(i\in \mathcal{I}')} 
\left| \widehat{E}\left(\sum_{i\in \mathcal{I}'}m_i \right)\right|^2=q^{|\mathcal{I}'|d}q^{-2d}|E|.$$
This implies that 
$$|R|\lesssim q^{\tfrac{kd}{2}}q^{-\tfrac{k(k-1)}{4}} {|E|}^{\tfrac{k+1}{2}}\sum_{\mathcal{I} \cup \mathcal{I}'= \{ 0, \dots ,k-1 \}} 
q^{-|\mathcal{I}'|(d+1)/2-|\mathcal{I}|} q^{|\mathcal{I}'|d}.$$

The largest term in the sum occurs when $\mathcal{I}=\emptyset$. We
conclude that
$$|R| \lesssim q^{\tfrac{kd}{2}} q^{-\tfrac{k(k+1)}{4}} {|E|}^
{\tfrac{k+1}{2}}.$$

The term $R$ is smaller than, say,  $\frac{M}{2}$ if
$$q^{\tfrac{kd}{2}} q^{-\tfrac{k(k+1)}{4}} {|E|}^{\tfrac{k+1}{2}}
\leq C {|E|}^{k+1}q^{-{k+1 \choose 2}},$$ with a sufficiently large
constant $C$, which happens if
$$ |E| \ge C'q^{\tfrac{k}{k+1}d}q^{\tfrac{k}{2}},$$ with a
sufficiently large constant $C'$ depending on the constants implicit
in the estimates above. This completes the proof.

\vskip.125in

\section{Proof of Lemma \ref{uptocrap}}

To prove Lemma \ref{uptocrap}, let $\pi_r(x)$ denote the $r$th coordinate of
$x$. There is no harm in assuming that $V_0=(0, \dots, 0)$. We may also
assume that $V_1, \dots, V_k$ are contained in ${\Bbb F}_q^k$. The condition
(\ref{equalnorm}) implies that
\begin{equation} \label{dotproduct} \sum_{r=1}^k \pi_r(V_i) \pi_r(V_j)=
\sum_{r=1}^k \pi_r(W_i) \pi_r(W_j). \end{equation}

Let $T$ be the linear transformation uniquely determined by the condition
$$ T(V_i)=V'_i.$$

In order to prove that $T$ is orthogonal, it suffices to show that
$$ ||Tx||=||x||$$ for any $x \not=(0, \dots, 0)$.

Since $V_j$s form a basis, by assumption, we have
$$ x=\sum_i t_i V_i, $$ so it suffices to show that
$$ ||x||=\sum_r \sum_{i,j} t_i t_j \pi_r(V_i) \pi_r(V_j)$$
$$=\sum_r \sum_{i,j} t_i t_j \pi_r(V'_i) \pi_r(V'_j)=||Tx||,$$ which follows
immediately from (\ref{dotproduct}).

Observe that we used the fact that orthogonality of $T$, the condition that
$T^t \cdot T=I$ is equivalent to the condition that $||Tx||=||x||$. To see
this observe that to show that $T^t \cdot T=I$ it suffices to show that
$T^tTx=x$ for all non-zero $x$. This, in turn, is equivalent to the
statement that
$$ <T^tTx,x>=||x||,$$ where
$$ <x,y>=\sum_{i=1}^k x_iy_i.$$

Now,
$$ <T^tTx,x>=<Tx, Tx>$$ by definition of the transpose, so the stated
equivalence is established. This completes the proof of Lemma
\ref{uptocrap}.

\vskip.125in

\section{Estimation of the Fourier transform of the sphere: proof of Lemma \ref{sphere}}
The proof of Lemma \ref{sphere} is fairly standard, but we outline
the argument for reader's convenience.  For any $m\in{\mathbb
F}^d_q$, we have
\begin{equation} \label{sphereparade}
\begin{array}{llllll} \widehat{S}_t(m)&=&
q^{-d} \sum_{x \in {\mathbb F}^d_q} q^{-1} \sum_{j \in {\mathbb F}_q} \chi(
j(\|x\|-t)) \chi( -  x \cdot m)\\ \hfill \\&=&q^{-1}\delta(m) +
q^{-d-1} \sum_{j \in {\mathbb F}^{*}_q} \chi(-jt) \sum_{x}
\chi( j\|x\|) \chi(- x \cdot m)\\ \hfill
\\&=&q^{-1}\delta(m)+   Q^d q^{-\frac{d+2}{2}} \sum_{j \in {\mathbb
F}^{*}_q}
\chi\left(\frac{\|m\|}{4j}+jt\right)\eta^d(-j),\end{array}\end{equation}
where the notation $\delta(m)=1$ if $m=(0\ldots,0)$ and $\delta(m)=0$
otherwise. 
In the last line we have completed the square, changed $j$ to
$-j$, and used $d$ times the Gauss sum equality
\begin{equation}
\sum_{c\in {\mathbb F}_q} \chi(jc^2) = \eta(j)\sum_{c\in{\mathbb
F}_q}\eta(c)\chi(c)=\eta(j)\sum_{c\in{\mathbb F}_q^*}\eta(c)\chi(c) =
Q\sqrt{q}\,\eta(j),
\label{gauss}\end{equation} where the constant $Q$ equals $\pm1$ or $\pm i$, depending on
$q$, and $\eta$ is the quadratic multiplicative character (or the
Legendre symbol) of ${\mathbb F}_q^*$.
The conclusion now follows from the following classical estimate due to
A. Weil (\cite{We48}).
\begin{theorem} \label{kloosterman} Let
$$ K(a)=\sum_{s \not=0} \chi(as+s^{-1}) \psi(s), $$ where, once
again, $\psi$ is a multiplicative character on ${\Bbb F}_q^{*}$. Then
$$ |K(a)| \leq 2 \sqrt{q}$$ if $a \not=0$.
\end{theorem}

\end{document}